\documentclass[11pt]{article}
\usepackage{amsmath,amsfonts}
\usepackage{theorem}
\input xy
\xyoption{all}
\usepackage{amstext}
\newcommand{\gC}{\Gamma}

\newcommand{\gc}{\gamma}

\newcommand{\gs}{\sigma}

\newcommand{\gp}{\psi}

\newcommand{\tN}{\tilde{N}}
\newcommand{\tgC}{\tilde{\gC}}

\newcommand{\p}{\prod}

\newcommand{\ra}{\rightarrow}
\newcommand{\Ra}{\Rightarrow}

\newcommand{\pf}{Proof: \ }
\newcommand{\epf}{\hfill$\spadesuit$\vspace{.25in}}

\newcommand{\Lim}[1]{\lim_{{#1} \ra \infty}}
\newcommand{\veC}[2]{({#1}_1,\ldots,{#1}_{{#2}})}

\newcommand{\Avr}[2]{\frac{1}{{#1}}\sum_{{#2}=1}^{{#1}}}

\newtheorem{lma}{Lemma}[section]
\newtheorem{cor}[lma]{Corollary}
\newtheorem{thm}[lma]{Theorem} 
 
\newtheorem{pro}[lma]{Proposition} 
 
\newtheorem{rmr}[lma]{Remark} 
\newtheorem{ntt}{Notation}

\begin{document}

\title{A Non Conventional Ergodic Theorem for a Nil-System} 
\author{T. Ziegler}
\date{}
\maketitle

\begin{abstract}  
We prove a non conventional pointwise convergence theorem for
a nilsystem, and give an explicit formula for the limit.
\end{abstract}

\section{Introduction}
Non conventional averages  of the form
\begin{equation}\label{e:MR}
  \Avr{N}{n}\p_{j=1}^kf_j( T^{jn}x)
\end{equation}
were first introduced by Furstenberg in the ergodic theoretic proof of 
Szemeredi's theorem on arithmetic progressions \cite{F}. 
Recent developments imply 
$L^2$ convergence of these averages to a limit for any $k$ \cite{HK}. 
We are interested
in pointwise convergence.
The case $k=1$ is the classical Birkhoff pointwise ergodic theorem.
J. Bourgain has proved pointwise convergence for the case $k=2$ \cite{Bo}. 
Pointwise convergence is 
not known for $k\ge 3$. 

Nilsystems arise in the study of the $L^2$ convergence of the averages 
(\ref{e:MR}). A nilspace is a homogeneous space of a nilpotent Lie 
group. A nilsystem consists of a finite measure nilspace $N/\gC$ with
a measure preserving transformation given by an element of $N$.  
It was shown by Furstenberg, Weiss \cite{FW}, Conze, Lesigne \cite{CL} that  
$2$-step nilsystems
`characterize' the $L^2$ behavior of the average (\ref{e:MR}) 
for the case $k=3$. By `characterize' we mean  that the $L^2$ average can 
be calculated by projecting on an appropriate factor - the characteristic 
factor of the average in question.
That $k$-step nil spaces play a role in the 
convergence of a $(k+1)$ length multiple recurrence average can be expected,
since for a $k$ step nilspace $N/\gC$, $a^{k+1}x\gC$ is a function of
$ax\gC,a^2x\gC,\ldots,a^{k}x\gC$, and this provides a constraint on the
$(k+1)$ tuple $T^nx,\ldots,T^{(k+1)}x$ in this case.
For this reason it is particularly interesting to study these averages
on nilsystems. 

\subsection{Statement of Result}
We prove a pointwise convergence theorem for nilsystem: 
Let $N$ be a $k$-step connected, simply connected nilpotent Lie group, 
$\gC$ a discrete subgroup s.t.
$N/\gC$ is compact. Let $N^1=N$, $N^i=[N^{i-1},N]$ ($N^{k+1}=1$), and
set  $\gC^i=\gC \cap N^i$. Then $\gC^i$ is a discrete subgroup
of $N^i$, and $N^i/\gC^i$ is compact (Malcev \cite{Ma}). 
Let $m_i$ be the probability measure
on $N^i/\gC^i$, invariant under translation by elements of $N^i$. Consider
the measure preserving system $(N/\gC,m_1,T)$, 
where $T$ denotes a translation by some element
$a \in N$.
In \cite{Le}, E. Lesigne proved the following theorem for step $2$
nilpotent groups:
\begin{thm}(Lesigne)
Suppose $T$ acts ergodically,
$f_1,f_2,f_3$ $\in L^{\infty}(N/\gC)$, 
then for almost all $x \in N$
\begin{equation*}
 \begin{split}
   \Lim{N}& \Avr{N}{n}T^nf_1(x\gC)T^{2n}f_{2}(x\gC)T^{3n}f_{3}(x\gC)  \\
    =& \int_{N/\gC}\int_{N^2/\gC^2}
       f_1(xy_1\gC)f_2(xy_1^2y_2\gC)
       f_{3}(xy_1^{3}y_2^{3}\gC)dm_1(y_1\gC^1)dm_2(y_2\gC^2)
 \end{split}
\end{equation*}
\end{thm} 
We generalize this result and prove the following ergodic theorem: 
\begin{thm}\label{main}
Suppose $T$ acts ergodically.
If $f_1, \ldots, f_{k+1}$ $\in L^{\infty}(N/\gC)$, 
then for almost all $x \in N$
\begin{equation*}
 \begin{aligned}
   \Lim{N}& \Avr{N}{n}\p_{j=1}^{k+1} T^{jn}f_j(x\gC) \\
   =& \int_{N/\gC} \ldots \int_{N^k/\gC^k}
       \p_{j=1}^{k+1} f_{j}(x\p_{i=1}^j y_i^{\binom{j}{i}}\gC) 
       \quad \p_{j=1}^{k+1} dm_j(y_j\gC^j).
 \end{aligned}
\end{equation*}
\end{thm}

The proof is a combination of the proof of Lesigne with the work
of Leibman \cite{LA} on geometric sequences in groups. For each $x$
we give the orbit of $(x,\ldots x)$ under $T \times \ldots \times T^{k+1}$
the structure of a nilsystem, with 
the action of a transformation $S_x$. The main technical
difficulty is with proving that the transformation $S_x$ 
is ergodic a.e. $x$.

\subsection{Acknowledgment}
I thank David Lehavi and Yoav Yaffe for helping me
make the proofs more readable.
This paper is a part of the author's PhD thesis. I thank 
my advisor Prof. Furstenberg for introducing me to ergodic theory, 
specifically to questions concerning non conventional ergodic theorems, 
and for fruitful discussions.  

\section{Proof of Theorem \ref{main}}
 
We define a set $\tN$ by  $\tN=N^1 \times \ldots \times N^k$.
For $\veC{x}{k}, \veC{y}{k} \in \tN$ we define
\begin{equation*}
 \begin{split}
   z_1&=x_1y_1, \\
   \p_{j=1}^{i}& z_j^{\binom{i}{j}}=\p_{j=1}^{i} x_j^{\binom{i}{j}}
     \p_{j=1}^{i} y_j^{\binom{i}{j}}.
 \end{split}
\end{equation*}
\[
  \veC{x}{k} \star \veC{y}{k} = \veC{z}{k},
\]
One should think of $\veC{x}{k}, \veC{y}{k}$ as representing the vectors 
\[
  (x_1,x_1^2 x_2, \ldots,\p_{j=1}^{k} x_j^{\binom{k}{j}}), \quad 
  (y_1,y_1^2 y_2, \ldots,\p_{j=1}^{k} y_j^{\binom{k}{j}})
\] 
respectively. Then
$\veC{z}{k}$ represents the coordinate product of these vectors.    
By Leibman \cite{LA}, $\tN$ is a group under this multiplication, thus
$\tN$ is a nilpotent connected simply connected Lie group.
Let $\tgC=\gC^1 \times \ldots \times \gC^k$, then $(\tgC,\star)$ is a discrete 
subgroup of $\tN$, and the quotient $\tN / \tgC$ is compact.
%
%
For $x \in N$ we define a transformation $S_x: \tN / \tgC \ra \tN / \tgC$ by:
\[
  S_x(\veC{y}{k}\star \tgC) 
     = ((a[a,x],e,\ldots,e) \star \veC{y}{k}) \star \tgC,
\]
and a function $I : \tN \ra (N / \gC)^{k+1}$ by:
\[
  I(\veC{y}{k}) 
  =(xy_1\gC,xy_1^2y_2\gC,\ldots,
    x\p_{j=1}^{k} y_j^{\binom{k+1}{j}}\gC).
\]
%
%
Let $\veC{\gc}{k} \in \tgC$ then we have for all $\veC{y}{k} \in \tN$:
\begin{equation*}
 \begin{split}
   I(\veC{y}{k} & \star \veC{\gc}{k})
      \\ 
   =& (xy_1\gc_1\gC,xy_1^2y_2\gc_1^2\gc_2\gC,\ldots,
       x\p_{j=1}^{k} y_j^{\binom{k+1}{j}}
        \p_{j=1}^{k} \gc_j^{\binom{k+1}{j}}\gC) 
      \\
   =& (xy_1\gC,xy_1^2y_2\gC,\ldots, x\p_{j=1}^{k} y_j^{\binom{k+1}{j}}\gC)
      \\ 
   =& I(\veC{y}{k}),
 \end{split}
\end{equation*}
thus $I$ can be defined on $\tN / \tgC$.
%
%
Let $\tilde{m}$ denote the probability measure on $\tN / \tgC$,
invariant under the action of $\tN$, then 
$\tilde{m}=m_1 \times \ldots \times m_k$.
Denote by $m_x$, $A_x$ the images of $\tilde{m}$, $\tN / \tgC$ under $I$. 
The function $I$ is an isomorphism of the systems
$(A_x,m_x, T \times T^2 \times \ldots \times T^{k+1})$ and
$(\tN / \tgC ,\tilde{m}, S_x)$ since
\begin{equation*}
 \begin{split}
   I(S_x & (\veC{y}{k} \star \tgC))
      \\
   =& I((a[a,x],e,\ldots,e) \star \veC{y}{k} \star \tgC)
      \\
   =& (xa[a,x]y_1\gC,x(a[a,x])^2y_1^2y_2\gC,\ldots,
       x(a[a,x])^{k+1}\p_{j=1}^{k} y_j^{\binom{k+1}{j}}\gC)
      \\
   =& (axy_1\gC,a^2xy_1^2y_2\gC,\ldots,
       a^{k+1}x\p_{j=1}^{k} y_j^{\binom{k+1}{j}}\gC)
      \\
   =& (T \times T^2 \times \ldots \times T^{k+1})
      (I(\veC{y}{k} \star \tgC )).
 \end{split}
\end{equation*}
If the action of $S_x$ is ergodic, then by Parry \cite{P1}, \cite{P2} 
the m.p.s 
$(\tN / \tgC ,\tilde{m}, S_x)$ is uniquely ergodic, and thus the m.p.s
$(A_x,m_x, T \times T^2 \times \ldots \times T^{k+1})$ is uniquely 
ergodic. For each such $x$ we have: for all $(y_1,\ldots,y_k)$ 
$\in \tN$, for all continuous functions $F$ on $(N / \gC)^{k+1}$,
\begin{equation*}
 \begin{split}
  \Lim{N}& \Avr{N}{n} F(a^nxy_1\gC,a^{2n}xy_1^2y_2\gC,\ldots,
                         a^{(k+1)n}x\p_{j=1}^{k} y_j^{\binom{k+1}{j}}\gC)\\
         =& \int_{N / \gC} \ldots \int_{N_k / \gC_k}
             F(xy_1\gC,xy_1^2y_2\gC,\ldots,
             x\p_{j=1}^{k} y_j^{\binom{k+1}{j}}\gC)          
             \quad \p_{j=1}^{k} dm_j(y_j\gC_j).                         
 \end{split}
\end{equation*}
In particular, for all such $F$
\begin{equation*}
 \begin{split}
  \Lim{N} \Avr{N}{n} & F(a^nx\gC,a^{2n}x\gC,\ldots,a^{(k+1)n}x\gC)   \\
                      =& \int_{N / \gC} \ldots \int_{N_k / \gC_k}
                         F(xy_1\gC,xy_1^2y_2\gC,\ldots,
                         x\p_{j=1}^{k} y_j^{\binom{k+1}{j}}\gC)       \\
                      \qquad& \qquad \qquad \qquad 
                       dm_1(y_1\gC)\ldots dm_k(y_k\gC).
 \end{split}
\end{equation*} 
If $S_x$ is ergodic then we proved the theorem for $f_1,\ldots,f_k$ continuous.
For passing to bounded measurable functions see \cite{CL} lemma 3 p. 173.

In order to finish the proof of the theorem, we must show that the
action of $S_x$ is ergodic for a.a. $x \in N$.
%
%
\subsection{Ergodicity of $S_x$}

We wish to show that the action of $S_x$ on $\tN$ is ergodic for
almost all (Haar) $x \in N$.
Let $\gs$ be a character of $\tN$ s.t. $\gs$ vanishes on $\tN' \star \tgC$.
By Green \cite{AGH} we must show that if $\gs \not \equiv 1$ then
\[
  m_h\{x \in N : \gs((a[a,x],e,\ldots,e))=1\} =0
\]
(as there is only countable number of characters on $\tN / \tN' \tgC$).
Suppose this set is of positive measure. Let $\gp(x)= 
\gs((a[a,x],e,\ldots,e))$, then $\gp$ is an analytic function
on $N$, therefore $\gp(x)=1$ on a set of positive measure implies 
$\gp \equiv 1$.
We define the functions $\gs_i$ on $N^i$ by
\[
  \gs_i(y_i)= \gs((e,\ldots,e,y_i,e,\ldots,e)).
\]
We have
\[
  \gs(\veC{y}{k})= \p_{i=1}^k \gs_i(y_i).
\]
For $y_1,\ldots,y_n \in N$ we define $R(y_1,\ldots,y_n)$ as follows:
\begin{equation*}
 \begin{split}
   R^0(y_1,\ldots,y_n) &= \{ y_1,\ldots,y_n \}                          \\
   R^j(y_1,\ldots,y_n) &= \bigcup_{i,l=0}^{j-1} \{ [R^{i},R^{l}] \} \\
     R(y_1,\ldots,y_n) &= \bigcup_{i=1}^{k} R^i(y_1,\ldots,y_n).
 \end{split}
\end{equation*}  
In other words, $R(y_1,\ldots,y_n)$ is the set of commutators involving
$y_1,\ldots,y_n$.
For all $x$
\begin{equation*}
 \begin{split}
  1=\gp(x)=\gs((a[a,x],e,\ldots,e))=\gs_1(a[a,x])                   \\
   =\gs_1(a)\gs_1([a,x])\gs_2(c_2(x)) \ldots \gs_k(c_k(x)), 
 \end{split}
\end{equation*}
where $c_i(x)$ are products of elements of $R(a,[a,x])$.
\begin{rmr}\label{multiplicity}
$\gs_i(xy) =\gs_i(x)\gs_i(y) $ modulo `corrections'
in $\gs_j$, $j>i$ computed on  commutators involving $x,y$.
\end{rmr}
Setting $x=e$ we get $\gs(a)=1$, thus
\begin{equation}\label{eq:[a,x]}
 \begin{split}
   \gs_1([a,x]) =&\gs_2(c_2(x))^{-1} \ldots \gs_k^{-1}(c_k(x))      \\
                =&\p_{i=2}^{k} \p_{m=1}^{l_i}\gs_i(r_{i_m}(x)), 
 \end{split}
\end{equation}
for all $x \in N$, where $r_{i_m}(x) \in R(a,[a,x])$ (as $R(a,[a,x])$ is 
closed under [,]).                                                  \\
We will show that $\gs_i$ are all trivial, and thus $\gs$ is trivial. 

\begin{ntt}
We use $[x_1,\ldots,x_n]$ for 
$[\ldots[[[x_1,x_2],x_3],x_4],\ldots,x_n]$.
\end{ntt}

\begin{ntt} 
Let $\cal{S}$ $= \{(\veC{s}{j}: 0 \le j \le k, 1 \le s_n \le k \}$.
For $A=\veC{s}{n} \in \cal{S}$ define $C_{A}$ by:
\begin{equation*}
  \begin{split}
     C_{A}^0&=a                                               \\
     C_{A}^j&= 
         \{[ C_{A}^{j-1},x_j]^{\pm 1}: x_j \in N^{s_j} \}     \\         
     C_{A}&= C_{A}^n
  \end{split}
\end{equation*}
\end{ntt}
{\em Example :} If $A=(1,3)$ then $C_A= \{ [[a,x]^{\pm 1},y]^{\pm 1} \}$
              where $x \in N$, $y \in N^3$.
\begin{ntt} 
 \begin{enumerate}
   \item For $A = \veC{s}{j} \in \cal{S}$, $|A| = \sum s_j$ 
   \item For $A \in \cal{S}$, $(i:A)$ means: $\forall v \in (C_A \cap N^i):
         \gs_i(v)=1$ (if $|A|<i-1$ the intersection is empty).
   \item $j^{\ge}$ ($>$) is some $l \ge j$ ($l >j$).
   \item Let $X$ be either $j$, $j^{\ge}$, $j^{>}$, then $X+$ is 
      a sequence of length $\ge 1$, whose entries are of type 
      $j$, $j^{\ge}$, $j^{>}$ respectively.
   \item For $x \in N^j$, $|x|=j$.
 \end{enumerate}
\end{ntt}


\begin{pro}\label{P:(i:A)}
 $\forall  A \in \cal{S}$, $1 \le i \le k$ : $(i:A)$.
\end{pro}
\pf
We will need the following lemma, which is based  on the fact that 
$\gs$ is multiplicative and thus $\gs \equiv 1$ on $\tN'$:


\begin{lma}\label{L:commutator}
Let $x \in N$, $y \in N^n$.
We define 
\[
  H(x,y)= \bigcup_{i=2}^{k} R^i(y_1,\ldots,y_n)  \quad
         (= R(x,y) \setminus  \{ [x,y]^{\pm 1} \}).
\]
If $\gs_i(h)=1$ $\forall i>1$, $h \in H(x,y) \cap N^i$,
then $\gs_m([x,y])^m=\gs_{m+1}([x,y])^{m+1}$, 
for all $m \le n$.
\end{lma}

\pf  
We denote by $\tilde{H}(x,y)$ the set of finite products of elements of
$H(x,y)$. By remark (\ref{multiplicity}) and since $H(x,y)$ is closed 
under commutators,
$\gs_i(\tilde{h})=1$ ($i>1$) for all $\tilde{h} \in \tilde{H}(x,y) \cap N^i$.
As $\gs|_{\tN'} \equiv 1$, we have
\begin{equation*}
 \begin{split}
  1=& \gs((x^{-1},e,\ldots,e)\star(e,\ldots,\stackrel{m}{y^{-1}},\ldots,e) 
      \star(x,e,\ldots,e)\star(e,\ldots,\stackrel{m}{y},\ldots,e))      \\
   =& \gs(z_1,\ldots,z_k) = \p_{i=1}^k \gs_i(z_i),
 \end{split}
\end{equation*}
we calculate the $z_i$: 
\begin{equation*}
  \begin{split}  
    z_1 &=\ldots=z_{m-1}=e,                                          \\
    z_m &=[x^m,y]=[x,y]^m \tilde{h_0}
 \end{split}
\end{equation*}  
for some $\tilde{h_0} \in \tilde{H}(x,y)$.
\begin{equation*}
  \begin{split}
     z_m^{(m+1)}z_{m+1}
      &=[x^{(m+1)},y^{(m+1)}]=[x,y]^{(m+1)^2} \tilde{h'_1} \qquad \Ra \\
     z_{m+1}&=[x,y]^{(m+1)} \tilde{h_1}
 \end{split}
\end{equation*}  
for some $\tilde{h_1},\tilde{h'_1} \in \tilde{H}(x,y)$.
We show by induction that $z_{m+j}=\tilde{h_j}$ for $j>1$,
$\tilde{h_j} \in \tilde{H}(x,y)$. For $j=2$:
\begin{equation*}
 \begin{split}
    z_m^{\binom{m+2}{2}}z_{m+1}^{(m+2)}z_{m+2}
     &=[x^{(m+2)},y^{\binom{m+2}{2}}]
     =[x,y]^{(m+1){\binom{m+2}{2}}}\tilde{h'_2} \qquad \Ra             \\
    z_{m+2}
     &=[x,y]^{-m{\binom{m+2}{2}}-(m+1)(m+2)+(m+2){\binom{m+2}{2}}}\tilde{h_2}
     =\tilde{h_2}
 \end{split}
\end{equation*}
for some $\tilde{h_2},\tilde{h'_2} \in \tilde{H}(x,y)$.
Suppose it is true for $2 \le i <j$, we have
\begin{equation*}
 \begin{split}
    z_m^{\binom{m+j}{j}}z_{m+1}^{\binom{m+j}{j-1}}\ldots z_{m+j}
     &=[x^{(m+j)},y^{\binom{m+j}{j}}]
     =[x,y]^{(m+j){\binom{m+j}{j}}}\tilde{h'_j} \qquad\Ra              \\
    z_{m+j}
     &=[x,y]^{-m{\binom{m+j}{j}}-(m+1){\binom{m+j}{j-1}}
      +(m+j){\binom{m+j}{j}}}\tilde{h_j}
     =\tilde{h_j}
 \end{split}
\end{equation*}
for some $\tilde{h_j},\tilde{h'_j} \in \tilde{H}(x,y)$.
Thus
\begin{equation*}
 \begin{split}
   1 =&\p_{i=1}^k \gs_i(z_i)
      =\gs_m([x,y]^m\tilde{h}_0)\gs_{m+1}([x,y]^{m+1}\tilde{h}_1)
       \p_{i=m+2}^k \gs_i(\tilde{h}_{i-m})                            \\
     =&\gs_m([x,y]^m)\gs_{m+1}([x,y]^{m+1}),   
 \end{split}
\end{equation*}
as $\tilde{H}([x,y],\tilde{h}_0)$, $\tilde{H}([x,y],\tilde{h}_1)$ 
$\subset \tilde{H}(x,y)$.
\epf

\begin{cor}
If $\gs_1([x,y])=1$, and $\gs_i(h)=1$ for all  $h \in H(x,y)\cap N^i$
then $\gs_i([x,y])=1$ for all $i$ s.t. $[x,y] \in N^i$.
\end{cor}
\pf
By lemma \ref{L:commutator}, $\gs_i([x,y])^i=1$. $N$ is connected, 
and  $\gs_i([e,y])=1$.
\epf



\begin{pro}\label{P:rules}
We prove the following rules:
 \begin{enumerate}
  \item $|A|\ge k \Longrightarrow \forall i: (i:A)$. (triviality rule)
  \item $\forall i: (i:j^{>} +) \Longrightarrow (1:j)$. (reducing rule)
  \item $(1:(A,l))$ \& $j > l$ \& $\forall i: (i:(A,l,j^{\ge} +))$ \&
        $\forall i: (i:(A,j-1,j^{\ge} +))$ 
        $\Longrightarrow$ $(1:(A,l,j-1))$ \&  $(1:(A,j-1,l))$.
        (expanding rule)
  \item $(1:(A,j))$ \&  $\forall i: (i:(A,j^{\ge},j^{\ge} +))$ \&
        $\forall i: (i:(A,j^{>}))$ $\Longrightarrow \forall i: (i:(A,j))$.
 \end{enumerate}
\end{pro}

\pf
Let $A=(s_1,\ldots,s_j)$, $0 \le j \le k$, $x_n \in N^{s_n}$.
 \begin{enumerate}
   \item $|A|\ge k$ implies that $[a,x_1,\ldots,x_j]=1$ thus 
         $\gs_i([a,x_1,\ldots,x_j])=1$.
   \item Let $x \in N^j$.                            \\
         We use equation \ref{eq:[a,x]}:
        \[
          \gs_1([a,x])= \p_{i=2}^{k} \p_{m=1}^{l_i}\gs_i(r_{i_m}(x)),
        \]
        where $r_{i_m}(x) \in R(a,[a,x])$.
        Denote $b=[a,x]$, then $r_{i_m}(x) \in R(a,b)$. 
        Looking at the $a$ appearing in a minimal number of
        brackets, $r_{i_m}(x)$ is of the form: $ [a,y_1,\dots, y_m]$
        where $m \ge 1$, $|y_n| \ge |b| >|x|$, and is thus of type 
        $(j^{> +})$. 
   \item Let $v=[a,x_1,\ldots,x_j]$, 
         $y \in N^l, z \in N^{j-1}$ then $yz \in N^l$.                      \\
         \begin{equation*}
          \begin{split}
            1 =& \gs_1([v,yz])=\gs_1([v,y][v,y,z][v,z])                 \\
              =& \gs_1([v,y])\gs_1([v,y,z])\gs_1([v,z])  
                 \p_{i=2}^k \p_{m=1}^{l_i} \gs_i(b_{i_m}(y,z,v))        \\
              =& \gs_1([v,y,z])\p_{i=2}^k \p_{m=1}^{l_i} \gs_i(b_{i_m}(y,z,v)),
          \end{split}
         \end{equation*}
         where $b_{i_m} \in R([v,y],[v,z],[v,y,z])$. 
         If $[v,y]$ appears in $b_{i_m}(y,z,v)$ then w.l.o.g.
         \[
         b_{i_m}(y,z,v)=[v,y,y_1,\ldots,y_m],
         \]
         where $m \ge 1$, $|y_n| \ge j$
         (look at the $[v,y]$ appearing in the minimal number of brackets)
         This is of type $(A,l,j^{\ge}+)$.
         Otherwise $[v,y]$ does not appear in $b_{i_m}(y,z,v)$, thus w.l.o.g. 
         \[
         b_{i_m}(y,z,v)=[v,z,z_1,\ldots,z_m],
         \]
         where $m \ge 1$, $|z_n| \ge j$. 
         This is of type $(A,j-1,j^{\ge}+)$.\\   
         For $(A,j-1,l)$, do the same with $[v,zy]$. 
   \item Let $v$ be of type $A$, $ x \in N^j$. By lemma \ref{L:commutator}, 
         we must show that for $h \in H(v,x)\cap N^i$, $\gs_i(h)=1$.
         But any such $h$ is either of the form $[v,r(v,x)]$ where 
         $r(v,x) \in R(v,x)$ and thus of type $(A,j^{>})$,  or 
         $[v,z_1,\ldots,z_m]$, $m \ge 1$, $|z_n| \ge j$, and thus of type 
         $(A,j^{\ge},j^{\ge} +)$.

\end{enumerate}


\begin{pro}\label{P:1toi}
$(1:(A,l)) \Longrightarrow (i:(A,l^{\ge} +))$, where $|A| \ge 0$, $l \ge 1$.
\end{pro}

\pf
We prove this by induction on $|A|$, $l$.
For $|A| = k$ any $l$ it is clear (as the commutator is trivial).
Suppose the statement is true for $|A| > d$, any $l$.
For $|A| = d$, $l=k-d$ it is clear (as the commutator is trivial).
Suppose the statement is true for $|A|=d$, any $j>l$. We follow the 
following scheme:

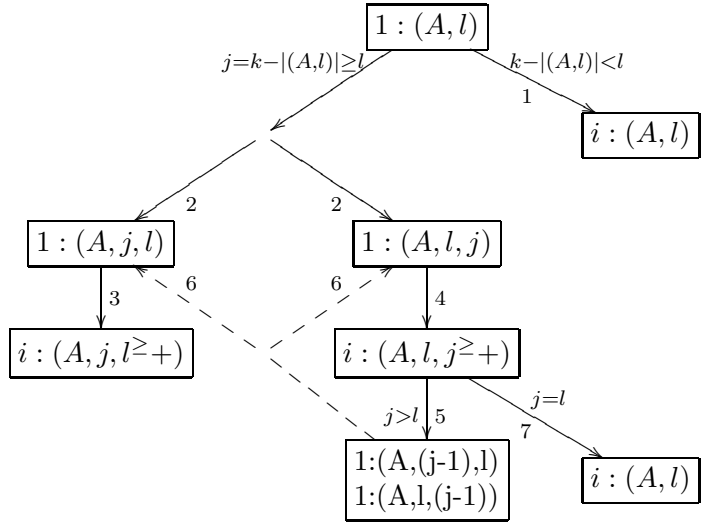
\xymatrix{
&  &*+[F]{ 1:(A,l) }  \ar[dl]_{j=k-|(A,l)| \ge l} \ar[dr]_1^{k-|(A,l)| < l}
\\
& \ar[dl]^2 \ar[dr]_2 & &  *+[F]{ i:(A,l) }                            \\
*+[F]{ 1:(A,j,l) } \ar[d]^3 & & *+[F]{ 1:(A,l,j) } \ar[d]^4    
\\
*+[F]{ i:(A,j,l^{\ge} +) } & \ar@{-->}[ul]_6  \ar@{-->}[ur]^6 & 
     *+[F]{ i:(A,l,j^{\ge} +) } \ar[d]^5_{j>l} \ar[dr]_7^{j=l}          \\
& & *+[F]{\txt{1:(A,(j-1),l) \\ 1:(A,l,(j-1))}} \ar@{--}[ul] & 
*+[F]{ i:(A,l) }
}
\ \\

Explanation for numbered arrows::
\begin{enumerate}
\item $k-|(A,l)| < l$ implies $(i:(A,l^{\ge},l^{\ge} +))$ by the 
      triviality rule $1$.
      $(1:(A,l))$ implies $(1:(A,j))$ for $j >l$ 
      (as $x \in N^j \Rightarrow x \in N^l$)
      By  the induction hypothesis $(i:(A,j))$ for $j >l$, and by rule 4
      $(i:(A,l))$.
\item By expanding rule 3 and triviality rule 1.
\item By the induction hypothesis (as $|A,j| >|A|$).
\item By the induction hypothesis (as $|A,l| >|A|$).
\item By expanding rule: $(j-1) \ge l$, $(1:(A,l))$, 
      $(i: (A,l,j^{\ge} +))$, and $(i:(A,(j-1),j^{\ge} +))$ as
      if $j-1=l$ it is the previous condition, and if $j-1 >l$ then
      $(1:(A,l))$ implies $(1:(A,j-1))$ which implies by the induction
      $(i:(A,j-1^{\ge} +))$.       
\item Repeat the procedure until $j=l$.
\item By this procedure we now know $(i:A,l^{\ge},l^{\ge} +)$.
      Repeat the argument in 1.

\end{enumerate}
\ \\
{\em Proof of proposition \ref{P:(i:A)}}
We prove this by induction on $j$: $(1:k)$ by triviality rule 1.
$(1:j) \Ra \forall i (i:j^{\ge} +)$ by proposition \ref{P:1toi} $(|A|=0)$.
By rule 2 this implies $(1:j-1)$.
\epf

\begin{pro}\label{P:sigma_i-equiv-1}
For all $i$, $\gs_i \equiv 1$.
\end{pro}

\pf
We prove this by induction on $k$.
By proposition \ref{P:(i:A)}, we have
\begin{equation}\label{eq:a-com}
 \gs_k([x_1,\ldots,\stackrel{m}{a},\ldots,x_{k-1}])=1
   \qquad \forall x_1,\ldots,x_{k-1} \in N, \quad 1 \le m \le k
\end{equation}  
Since $\gs|_{\tgC} \equiv 1$,
\[ 
  \gs_k([\gc_1,\ldots,\gc_k])=1 \qquad \forall \gc_1,\ldots,\gc_k \in \gC.
\]
From (\ref{eq:a-com}) 
\[
  \gs_k([a,\gc_2,\ldots,\gc_k])=1 \qquad \forall \gc_2,\ldots,\gc_k \in \gC.
\]
Since $\gs_k$ is a character on $N^k$, and $a$ acts ergodically on $N/\gC$,
we have
\[
  \gs_k([x_1,\gc_2,\ldots,\gc_k])=1 \qquad 
  \forall \gc_2,\ldots,\gc_k \in \gC, \quad x_1 \in N.
\]
Using (\ref{eq:a-com}) and induction, we get
\[
 \gs_k([x_1,\ldots,x_k])=1 \qquad \forall x_1,\ldots,x_k \in N,
\]
thus $\gs_k \equiv 1$.
Now $\gs_k,\ldots,\gs_{j+1} \equiv 1$ implies that $\gs_j$ is a character.
For $j>1$, by Proposition \ref{P:(i:A)}, we have
\begin{equation*}
 \gs_j([x_1,\ldots,\stackrel{m}{a},\ldots,x_{j-1}])=1
   \qquad \forall x_1,\ldots,x_{j-1} \in N, \quad 1 \le m \le j ,
\end{equation*}
and by the same argument as above we have $\gs_j \equiv 1$.
For $j=1$, $\gs_1$ is a character which satisfies
\[
  \gs_1(a)=1 , \qquad \gs_1|_\gC \equiv 1,
\]
thus $\gs_1 \equiv 1$.   
\epf



\begin{thebibliography}{99}
\bibitem{AGH} Auslander, L., Green, L., Hahn, F.  
        {\em Flows on homogeneous spaces}.
        Annals of Mathematics Studies, No. $53$ (1963).
\bibitem{Bo} Bourgain, J.  
        {\em Pointwise ergodic theorems for arithmetic sets}. 
        Inst. Hautes Études Sci. Publ. Math. No. $69$: $5$-$45$ (1989).
\bibitem{CL} Conze, J.P., Lesigne, E. 
        {\em Th{\'e}or{\`e}mes ergodiques pour des mesures diagonales}.
        Bull. Soc. Math. France $112$ no. $2$: $143$-$175$ (1984).
\bibitem{F} Furstenberg, H.  
        {\em Ergodic behavior of diagonal measures and a theorem of 
         Szemeredi on arithmetic progressions}. 
        J. Analyse Math. 31: $204$-$256$ (1977). 
\bibitem{FW} Furstenberg, H., Weiss B. 
        {\em A mean ergodic theorem for 
         $(1/N)\sum^N_{n=1}f (T^n x)g(T^ {n\sp 2}x)$}. 
        Convergence in ergodic theory and probability 
        Ohio State Univ. Math. Res. Inst. Publ., 5, de Gruyter, Berlin:
        $193$-$227$ (1996).
\bibitem{HK} Host, B. , Kra, B : personal communication.
\bibitem{LA} Leibman, A 
        {\em Polynomial sequences in groups}.
        Journal of Algebra, $201$: $189$-$206$ (1998). 
\bibitem{Le} Lesigne, E.  
        {\em Theorems ergodiques pour une translation sur une nilvariete}
        Erg. Thm and Dyn.Sys, $9(1)$: $115$-$126$ (1989).
\bibitem{Ma} Malcev, A.I.
        {\em On a class of homogeneous spaces}.
        Amer. Math. Soc. Translation no. $39$ (1951).
\bibitem{P1} Parry, W. (1969)
        {\em Ergodic properties of affine transformations and flows on 
         nilmanifolds}. 
        Amer. J. Math. $91$ :$757$-$771$.
\bibitem{P2}  Parry, W.
        {\em Dynamical systems on nilmanifolds}.
        Bull. London Math. Soc. $2$: $37$-$40$ (1970).


\end{thebibliography}
\end{document}